\documentclass[11pt]{amsart}
\usepackage{amssymb,amsmath}
\usepackage{graphicx}
\usepackage[all]{xy}

\title[Adding new T-homotopy equivalences]{T-homotopy and refinement of observation (II) : Adding new T-homotopy equivalences}
\author[P. Gaucher]{Philippe Gaucher}
\address{Preuves Programmes et Syst{\`e}mes\\ Universit{\'e} Paris 7--Denis Diderot\\
Case 7014\\2 Place Jussieu\\ 75251 PARIS Cedex 05\\ France}
\email{gaucher@pps.jussieu.fr}
\urladdr{http://www.pps.jussieu.fr/{\~{}}gaucher/}

\newcommand{\C}{\mathcal{C}}

\newcommand{\p}\times
\renewcommand{\vec}{\overrightarrow}
\renewcommand{\P}{\mathbb{P}}

\newcommand{\be}{\begin{equation}}
\newcommand{\ee}{\end{equation}}
\newcommand{\bea}{\begin{eqnarray}}
\newcommand{\eea}{\end{eqnarray}}
\newcommand{\beas}{\begin{eqnarray*}}
\newcommand{\eeas}{\end{eqnarray*}}

\newtheorem*{thmN}{Theorem}

\newtheorem{thm}{Theorem}[section]
\newtheorem{prop}[thm]{Proposition}
\newtheorem{lem}[thm]{Lemma}

\newtheorem{defn}[thm]{Definition}

\newtheorem{nota}[thm]{Notation}
\newcommand{\bd}{\begin{defn}}
\newcommand{\ed}{\end{defn}}
\newcommand{\bp}{\begin{prop}}
\newcommand{\ep}{\end{prop}}
\newcommand{\bth}{\begin{thm}}
\renewcommand{\eth}{\end{thm}}
\newcommand{\bpf}{\begin{proof}}
\newcommand{\epf}{\end{proof}}

\newcommand{\fl}[1]{\ar@{->}[l]_{#1}}
\newcommand{\fr}[1]{\ar@{->}[r]^-{#1}}
\newcommand{\fd}[1]{\ar@{->}[d]_{#1}}
\newcommand{\fu}[1]{\ar@{->}[u]^{#1}}
\newcommand{\f}[2]{\ar@{->}[#1]|{#2}}
\newcommand{\ff}[2]{\ar@2{->}[#1]|{#2}}
\newcommand{\frr}[1]{\ar@{->}[rr]^{#1}}

\renewcommand{\top}{{\mathbf{Top}}}
\newcommand{\poset}{{\brm{PoSet}}}

\newcommand{\iso}{\cong}
\newcommand{\lp}{\left(}
\newcommand{\rp}{\right)}

\newcommand{\vI}{\vec{I}}

\renewcommand{\leq}{\leqslant}
\renewcommand{\geq}{\geqslant}

\def\cartesien{%
  \ar@{-}[]+R+<6pt,-2pt>;[]+RD+<6pt,-6pt>%
  \ar@{-}[]+D+<2pt,-6pt>;[]+RD+<6pt,-6pt>%
}
\def\cocartesien{%
  \ar@{-}[]+L+<-6pt,+2pt>;[]+LU+<-6pt,+6pt>%
  \ar@{-}[]+U+<-2pt,+6pt>;[]+LU+<-6pt,+6pt>%
}

\newcommand{\brm}[1]{\rm{\mathbf{#1}}}
\renewcommand{\top}{{\brm{Top}}}
\newcommand{\gltop}{{\brm{glTop}}}
\newcommand{\dtop}{{\brm{Flow}}}
\newcommand{\set}{{\brm{Set}}}
\newcommand{\tgltop}{{\brm{glTOP}}}
\newcommand{\glob}{{\rm{Glob}}}

\newcommand{\liminj}{\varinjlim}

\DeclareMathOperator{\id}{Id}

\DeclareMathOperator{\cell}{{\brm{cell}}}
\DeclareMathOperator{\cof}{{\brm{cof}}}
\DeclareMathOperator{\inj}{{\brm{inj}}}
\newcommand{\hda}{{\cell(\dtop)}}

\begin{document}

\begin{abstract} 
  This paper is the second part of a series of papers about a new
  notion of T-homotopy of flows. It is proved that the old definition
  of T-homotopy equivalence does not allow the identification of the
  directed segment with the $3$-dimensional cube. This contradicts a
  paradigm of dihomotopy theory. A new definition of T-homotopy
  equivalence is proposed, following the intuition of refinement of
  observation. And it is proved that up to weak S-homotopy, an old
  T-homotopy equivalence is a new T-homotopy equivalence. The
  left-properness of the weak S-homotopy model category of flows is
  also established in this second part.  The latter fact is used
  several times in the next papers of this series.
\end{abstract}

\maketitle

\tableofcontents

\section{Outline of the paper}

\begin{figure}
\begin{center}
\includegraphics[width=5cm]{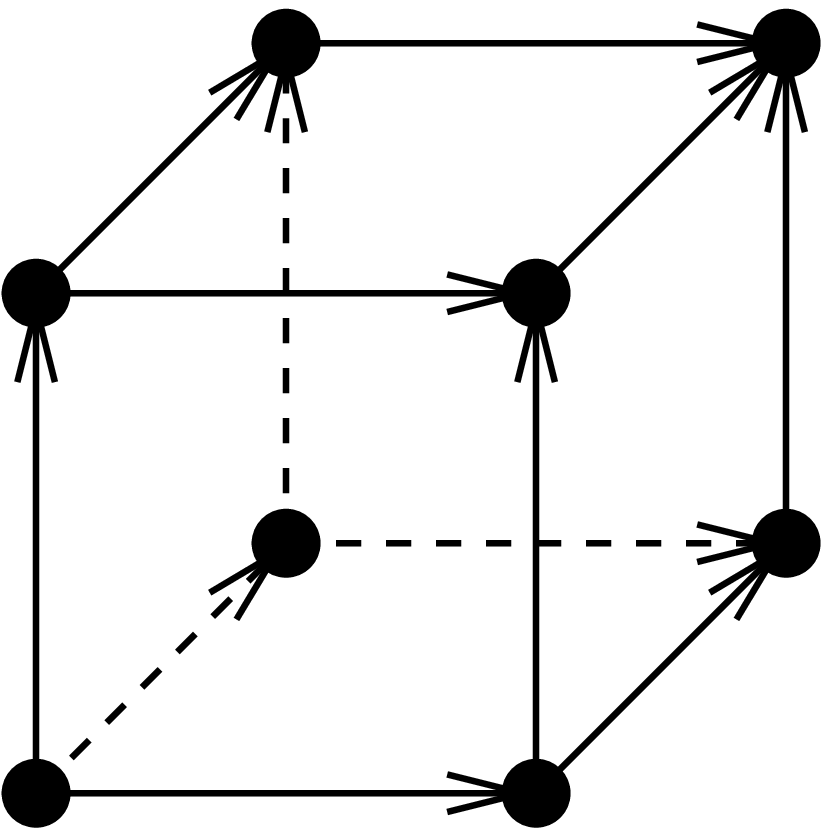}
\end{center}
\caption{The full $3$-cube}
\label{3cube}
\end{figure}

The first part \cite{1eme} of this series was an expository paper
about the geometric intuition underlying the notion of T-homotopy. The
purpose of this second paper is to prove that the class of old
T-homotopy equivalences introduced in \cite{diCW} and in \cite{model2}
is actually not big enough. Indeed, the only kind of old T-homotopy
equivalence consists of the deformations which locally act like in
Figure~\ref{ex1}. So it becomes impossible with this old definition to
identify the directed segment of Figure~\ref{ex1} with the full
$3$-cube of Figure~\ref{3cube} by a zig-zag sequence of weak
S-homotopy and of T-homotopy equivalences preserving the initial state
and the final state of the $3$-cube since every point of the $3$-cube
is related to three distinct edges. This contradicts the fact that
concurrent execution paths cannot be distinguished by observation.
The end of the paper proposes a new definition of T-homotopy
equivalence following the paradigm of invariance by refinement of
observation. It will be checked that the preceding drawback is then
overcome.

This second part gives only a motivation for the new definition of
T-homotopy. Further developments and applications are given in
\cite{3eme}, \cite{4eme} and \cite{hocont}. The left-properness of the
model category structure of \cite{model3} is also established in this
paper. The latter result is used several times in the next papers of
this series (e.g., \cite{3eme} Theorem~11.2, \cite{4eme} Theorem~9.2).

\begin{figure}
\[
\xymatrix{\widehat{0} \ar@{->}[rrrr]^-{U} &&&& \widehat{1} \\
\widehat{0} \ar@{->}[rr]^-{U'} && A \ar@{->}[rr]^-{U''} && \widehat{1}}
\]
\caption{The simplest example of T-homotopy equivalence}
\label{ex1}
\end{figure}

Section~\ref{gldef} collects some facts about globular complexes and
their relationship with the category of flows. Indeed, it is not known
how to establish the limitations of the old form of T-homotopy
equivalence without using globular complexes together with a
compactness argument.  Section~\ref{rappelT} recalls the notion of old
T-homotopy equivalence of flows which is a kind of morphism between
flows coming from globular complexes (the class of flows $\hda$).
Section~\ref{factI} presents elementary facts about relative
$I^{gl}_+$-cell complexes which will be used later in the paper.
Section~\ref{sectionleftproperflow} proves that the model category of
flows is left proper. This technical fact is used in the proof of the
main theorem of the paper, and it was not established in
\cite{model3}. Section~\ref{main1} proves the first main theorem of
the paper:
\begin{thmN} (Theorem~\ref{pasiso}) Let $n\geq 3$. There does not
  exist any zig-zag sequence of S-homotopy equivalences and of old
  T-homotopy equivalences between the flow associated with the
  $n$-cube and the flow associated with the directed segment.
\end{thmN} Finally Section~\ref{main2} proposes a new definition of
T-homotopy equivalence and the second main theorem of the paper is
proved:
\begin{thmN} (Theorem~\ref{OKnew}) Every T-homotopy in the old sense
  is the composite of an S-homotopy equivalence with a T-homotopy
  equivalence in the new sense~\footnote{Since a T-homotopy in the old
    sense is a T-homotopy in the new sense only up to S-homotopy, the
    terminology ``generalized T-homotopy'' used in Section~\ref{main2}
    may not be the best one. However, this terminology is used in the
    other papers of this series so we keep it to avoid any
    confusion.}.
\end{thmN}

\section{Prerequisites and notations}

The initial object (resp. the terminal object) of a category $\C$, if
it exists, is denoted by $\varnothing$ (resp. $\mathbf{1}$).

Let $\C$ be a cocomplete category.  If $K$ is a set of morphisms of
$\C$, then the class of morphisms of $\C$ that satisfy the RLP
(\textit{right lifting property}) with respect to any morphism of $K$
is denoted by $\inj(K)$ and the class of morphisms of $\C$ that are
transfinite compositions of pushouts of elements of $K$ is denoted by
$\cell(K)$. Denote by $\cof(K)$ the class of morphisms of $\C$ that
satisfy the LLP (\textit{left lifting property}) with respect to the
morphisms of $\inj(K)$.  It is a purely categorical fact that
$\cell(K)\subset \cof(K)$. Moreover, every morphism of $\cof(K)$ is a
retract of a morphism of $\cell(K)$ as soon as the domains of $K$ are
small relative to $\cell(K)$ (\cite{MR99h:55031} Corollary~2.1.15). An
element of $\cell(K)$ is called a \textit{relative $K$-cell complex}.
If $X$ is an object of $\C$, and if the canonical morphism
$\varnothing\longrightarrow X$ is a relative $K$-cell complex, then
the object $X$ is called a \textit{$K$-cell complex}.

Let $\C$ be a cocomplete category with a distinguished set of
morphisms $I$. Then let $\cell(\C,I)$ be the full subcategory of $\C$
consisting of the object $X$ of $\C$ such that the canonical morphism
$\varnothing\longrightarrow X$ is an object of $\cell(I)$. In other
terms, $\cell(\C,I)=(\varnothing\!\downarrow \! \C) \cap \cell(I)$.

It is obviously impossible to read this paper without a strong
familiarity with \textit{model categories}. Possible references for
model categories are \cite{MR99h:55031}, \cite{ref_model2} and
\cite{MR1361887}.  The original reference is \cite{MR36:6480} but
Quillen's axiomatization is not used in this paper. The axiomatization
from Hovey's book is preferred.  If $\mathcal{M}$ is a
\textit{cofibrantly generated} model category with set of generating
cofibrations $I$, let $\cell(\mathcal{M}) := \cell(\mathcal{M},I)$:
this is the full subcategory of \textit{cell complexes} of the model
category $\mathcal{M}$. A cofibrantly generated model structure
$\mathcal{M}$ comes with a \textit{cofibrant replacement functor}
$Q:\mathcal{M} \longrightarrow \cell(\mathcal{M})$. In all usual model
categories which are cellular (\cite{ref_model2} Definition~12.1.1),
all the cofibrations are monomorphisms.  Then for every monomorphism
$f$ of such a model category $\mathcal{M}$, the morphism $Q(f)$ is a
cofibration, and even an inclusion of subcomplexes (\cite{ref_model2}
Definition~10.6.7) because the cofibrant replacement functor $Q$ is
obtained by the small object argument, starting from the identity of
the initial object. This is still true in the model category of flows
reminded in Section~\ref{remflow} since the class of cofibrations
which are monomorphisms is closed under pushout and transfinite
composition.

A \textit{partially ordered set} $(P,\leq)$ (or \textit{poset}) is a
set equipped with a reflexive antisymmetric and transitive binary
relation $\leq$. A poset is \textit{locally finite} if for any
$(x,y)\in P\p P$, the set $[x,y]=\{z\in P,x\leq z\leq y\}$ is finite.
A poset $(P,\leq)$ is \textit{bounded} if there exist $\widehat{0}\in
P$ and $\widehat{1}\in P$ such that $P = [\widehat{0},\widehat{1}]$
and such that $\widehat{0} \neq \widehat{1}$. For a bounded poset $P$,
let $\widehat{0}=\min P$ (the bottom element) and $\widehat{1}=\max P$
(the top element). In a poset $P$, the interval $]\alpha,-]$ (the
sub-poset of elements of $P$ strictly bigger than $\alpha$) can also
be denoted by $P_{>\alpha}$.

A poset $P$, and in particular an ordinal, can be viewed as a small
category denoted in the same way: the objects are the elements of $P$
and there exists a morphism from $x$ to $y$ if and only if $x\leq y$.
If $\lambda$ is an ordinal, a \textit{$\lambda$-sequence} in a
cocomplete category $\C$ is a colimit-preserving functor $X$ from
$\lambda$ to $\C$. We denote by $X_\lambda$ the colimit $\liminj X$
and the morphism $X_0\longrightarrow X_\lambda$ is called the
\textit{transfinite composition} of the morphisms
$X_\mu\longrightarrow X_{\mu+1}$.

A model category is \textit{left proper} if the pushout of a weak
equivalence along a cofibration is a weak equivalence. The model
categories $\top$ and $\dtop$ (see below) are both left proper (cf.
Theorem~\ref{leftproperflow} for $\dtop$).

In this paper, the notation $\xymatrix@1{\ar@{^{(}->}[r]&}$ means
\textit{cofibration}, the notation $\xymatrix@1{\ar@{->>}[r]&}$ means
\textit{fibration}, the notation $\simeq$ means \textit{weak equivalence}, 
and the notation $\iso$ means \textit{isomorphism}.

\section{Reminder about the category of flows}
\label{remflow}

The category $\top$ of \textit{compactly generated topological spaces}
(i.e. of weak Hausdorff $k$-spaces) is complete, cocomplete and
cartesian closed (more details for this kind of topological spaces in
\cite{MR90k:54001,MR2000h:55002}, the appendix of \cite{Ref_wH} and
also the preliminaries of \cite{model3}). For the sequel, any
topological space will be supposed to be compactly generated. A
\textit{compact space} is always Hausdorff.

The category $\top$ is equipped with the unique model structure having
the \textit{weak homotopy equivalences} as weak equivalences and
having the \textit{Serre fibrations}~\footnote{that is a continuous
map having the RLP with respect to the inclusion $\mathbf{D}^n\p
0\subset \mathbf{D}^n\p [0,1]$ for any $n\geq 0$ where $\mathbf{D}^n$
is the $n$-dimensional disk.} as fibrations.

The time flow of a higher dimensional automaton is encoded in an
object called a \textit{flow} \cite{model3}. A flow $X$ contains a set
$X^0$ called the \textit{$0$-skeleton} whose elements correspond to
the \textit{states} (or \textit{constant execution paths}) of the
higher dimensional automaton. For each pair of states
$(\alpha,\beta)\in X^0\p X^0$, there is a topological space
$\P_{\alpha,\beta}X$ whose elements correspond to the
\textit{(non-constant) execution paths} of the higher dimensional
automaton \textit{beginning at} $\alpha$ and \textit{ending at}
$\beta$. For $x\in \P_{\alpha,\beta}X$, let $\alpha=s(x)$ and
$\beta=t(x)$. For each triple $(\alpha,\beta,\gamma)\in X^0\p X^0\p
X^0$, there exists a continuous map $*:\P_{\alpha,\beta}X\p
\P_{\beta,\gamma}X\longrightarrow \P_{\alpha,\gamma}X$ called the
\textit{composition law} which is supposed to be associative in an
obvious sense. The topological space $\P
X=\bigsqcup_{(\alpha,\beta)\in X^0\p X^0}\P_{\alpha,\beta}X$ is called
the \textit{path space} of $X$. The category of flows is denoted by
$\dtop$. A point $\alpha$ of $X^0$ such that there are no non-constant
execution paths ending at $\alpha$ (resp. starting from $\alpha$) is
called an \textit{initial state} (resp. a \textit{final state}). A
morphism of flows $f$ from $X$ to $Y$ consists of a set map $f^0:X^0
\longrightarrow Y^0$ and a continuous map $\P f: \P X \longrightarrow
\P Y$ preserving the structure. A flow is therefore ``almost'' a small
category enriched in $\top$. A flow $X$ is \textit{loopless} if for
every $\alpha\in X^0$, the space $\P_{\alpha,\alpha}X$ is empty.

\begin{figure}
\begin{center}
\includegraphics[width=7cm]{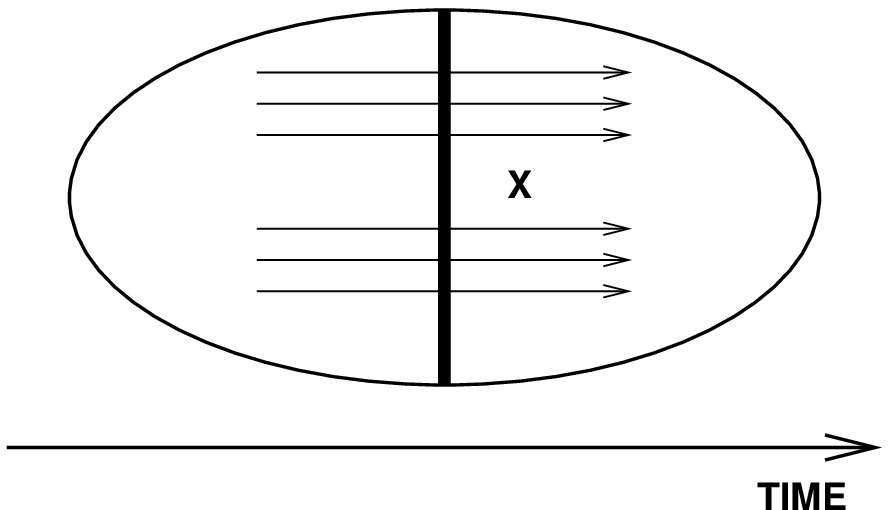}
\end{center}
\caption{Symbolic representation of
$\glob(Z)$ for some topological space $Z$} \label{exglob}
\end{figure}

Here are four fundamental examples of flows:
\begin{enumerate}
\item Let $S$ be a set. The flow associated with $S$, still denoted by
  $S$, has $S$ as set of states and the empty space as path space.
  This construction induces a functor $\set \rightarrow \dtop$ from
  the category of sets to that of flows. The flow associated with a
  set is loopless. 
\item Let $(P,\leq)$ be a poset. The flow associated with $(P,\leq)$,
  and still denoted by $P$ is defined as follows: the set of states of
  $P$ is the underlying set of $P$; the space of morphisms from
  $\alpha$ to $\beta$ is empty if $\alpha\geq \beta$ and equals to
  $\{(\alpha,\beta)\}$ if $\alpha<\beta$ and the composition law is
  defined by $(\alpha,\beta)*(\beta,\gamma) = (\alpha,\gamma)$. This
  construction induces a functor $\poset \rightarrow \dtop$ from the
  category of posets together with the strictly increasing maps to the
  category of flows. The flow associated with a poset is loopless. 
\item The flow $\glob(Z)$ defined by \beas
  && \glob(Z)^0=\{\widehat{0},\widehat{1}\} \\
  && \P \glob(Z)=Z \hbox{ with }s(z)=\widehat{0}\hbox{ and }t(z)=\widehat{1}\hbox{ for all }z\in Z 
\eeas and a trivial composition law (cf.
  Figure~\ref{exglob}). It is called the \textit{globe} of $Z$. 
\item The \textit{directed segment} $\vI$ is by definition
  $\glob(\{0\}) \iso \{\widehat{0} < \widehat{1}\}$. 
\end{enumerate}

The category $\dtop$ is equipped with the unique model structure
such that \cite{model3}: 
\begin{itemize}
\item The weak equivalences are the \textit{weak S-homotopy equivalences}, 
i.e. the morphisms of flows $f:X\longrightarrow Y$ such that
$f^0:X^0\longrightarrow Y^0$ is a bijection and such that $\P f:\P
X\longrightarrow \P Y$ is a weak homotopy equivalence. 
\item The fibrations are the morphisms of flows
$f:X\longrightarrow Y$ such that $\P f:\P X\longrightarrow \P Y$ is a
Serre fibration.
\end{itemize}
This model structure is cofibrantly generated. The set of generating
cofibrations is the set $I^{gl}_+=I^{gl}\cup
\{R:\{0,1\}\longrightarrow \{0\},C:\varnothing\longrightarrow \{0\}\}$
with
\[I^{gl}=\{\glob(\mathbf{S}^{n-1})\subset \glob(\mathbf{D}^{n}), n\geq
0\}\] where $\mathbf{D}^{n}$ is the $n$-dimensional disk and
$\mathbf{S}^{n-1}$ the $(n-1)$-dimensional sphere. The set of
generating trivial cofibrations is
\[J^{gl}=\{\glob(\mathbf{D}^{n}\p\{0\})\subset
\glob(\mathbf{D}^{n}\p [0,1]), n\geq 0\}.\]

If $X$ is an object of $\cell(\dtop)$, then a presentation of the
morphism $\varnothing \longrightarrow X$ as a transfinite composition
of pushouts of morphisms of $I^{gl}_+$ is called a \textit{globular
decomposition} of $X$.

\section{Globular complex}
\label{gldef}

The reference is \cite{model2}. A \textit{globular complex} is a
topological space together with a structure describing the sequential
process of attaching \textit{globular cells}. A general globular
complex may require an arbitrary long transfinite construction. We
restrict our attention in this paper to globular complexes whose
globular cells are morphisms of the form
$\glob^{top}(\mathbf{S}^{n-1}) \longrightarrow
\glob^{top}(\mathbf{D}^{n})$ (cf. Definition~\ref{globe_top}).

\bd 
A {\rm multipointed topological space} $(X,X^0)$ is a pair of
topological spaces such that $X^0$ is a discrete subspace of $X$.  A
morphism of multipointed topological spaces $f:(X,X^0)\longrightarrow
(Y,Y^0)$ is a continuous map $f:X\longrightarrow Y$ such that
$f(X^0)\subset Y^0$. The corresponding category is denoted by
$\top^m$. The set $X^0$ is called the {\rm $0$-skeleton} of $(X,X^0)$.
The space $X$ is called the {\rm underlying topological space} of
$(X,X^0)$. 
\ed

The category of multipointed spaces is cocomplete.

\bd \label{globe_top}
Let $Z$ be a topological space. The {\rm globe of $Z$}, which is
denoted by $\glob^{top}(Z)$, is the multipointed space
\[(|\glob^{top}(Z)|,\{\widehat{0},\widehat{1}\})\] where the topological space
$|\glob^{top}(Z)|$ is the quotient of
$\{\widehat{0},\widehat{1}\}\sqcup (Z\p[0,1])$ by the relations
$(z,0)=(z',0)=\widehat{0}$ and $(z,1)=(z',1)=\widehat{1}$ for any
$z,z'\in Z$. In particular, $\glob^{top}(\varnothing)$ is the
multipointed space
$(\{\widehat{0},\widehat{1}\},\{\widehat{0},\widehat{1}\})$.
\ed

\begin{nota} 
  If $Z$ is a singleton, then the globe of $Z$ is denoted by
  $\vI^{top}$.
\end{nota}

\bd 
Let $I^{gl,top}:=\{\glob^{top}(\mathbf{S}^{n-1})\longrightarrow
\glob^{top}(\mathbf{D}^{n}),n\geq 0\}$. 
A {\rm relative globular precomplex} is a relative $I^{gl,top}$-cell
complex in the category of multipointed topological spaces. 
\ed

\bd 
A {\rm globular precomplex} is a $\lambda$-sequence of multipointed
topological spaces $X:\lambda\longrightarrow \top^m$ such that $X$ is
a relative globular precomplex and such that $X_0=(X^0,X^0)$ with
$X^0$ a discrete space. This $\lambda$-sequence is characterized by a
presentation ordinal $\lambda$, and for any $\beta<\lambda$, an
integer $n_\beta\geq 0$ and an attaching map
$\phi_\beta:\glob^{top}(\mathbf{S}^{n_\beta-1}) \longrightarrow
X_\beta$. The family $(n_\beta,\phi_\beta)_{\beta<\lambda}$ is called
the {\rm globular decomposition} of $X$.
\ed

Let $X$ be a globular precomplex. The $0$-skeleton of $\liminj X$ is
equal to $X^0$.

\bd A morphism of globular precomplexes $f:X\longrightarrow Y$ is a
morphism of multipointed spaces still denoted by $f$ from $\liminj X$
to $\liminj Y$. \ed

\begin{nota} 
If $X$ is a globular precomplex, then the underlying topological space
of the multipointed space $\liminj X$ is denoted by $|X|$ and the
$0$-skeleton of the multipointed space $\liminj X$ is denoted by
$X^0$. 
\end{nota}

\bd 
Let $X$ be a globular precomplex. The space $|X|$ is called the {\rm
underlying topological space} of $X$. The set $X^0$ is called the {\rm
$0$-skeleton} of $X$. 
\ed

\bd 
Let $X$ be a globular precomplex. A morphism of globular precomplexes
$\gamma:\vI^{top}\longrightarrow X$ is a {\rm non-constant execution
path} of $X$ if there exists $t_0=0<t_1<\dots<t_{n}=1$ such that:
\begin{enumerate}
\item  $\gamma(t_i)\in X^0$ for any $0 \leq i \leq n$, 
\item $\gamma(]t_i,t_{i+1}[)\subset \glob^{top}(\mathbf{D}^{n_{\beta_i}} 
\backslash \mathbf{S}^{n_{\beta_i}-1})$ for some $(n_{\beta_i},\phi_{\beta_i})$ 
of the globular decomposition of $X$, 
\item  for $0\leq i<n$, there exists $z^i_\gamma\in \mathbf{D}^{n_{\beta_i}}\backslash 
\mathbf{S}^{n_{\beta_i}-1}$ and a strictly increasing continuous map
$\psi^i_\gamma:[t_i,t_{i+1}]\longrightarrow [0,1]$ such that
$\psi^i_\gamma(t_i)=0$ and $\psi^i_\gamma(t_{i+1})=1$ and for any
$t\in [t_i,t_{i+1}]$, $\gamma(t)=(z^i_\gamma,\psi^i_\gamma(t))$.
\end{enumerate}
In particular, the restriction $\gamma\!\restriction_{]t_i,t_{i+1}[}$
of $\gamma$ to $]t_i,t_{i+1}[$ is one-to-one. The set of non-constant
execution paths of $X$ is denoted by ${\P}^{top}(X)$. 
\ed

\bd 
A morphism of globular precomplexes $f:X\longrightarrow Y$ is {\rm
non-decreasing} if the canonical set map
$\top([0,1],|X|)\longrightarrow \top([0,1],|Y|)$ induced by
composition by $f$ yields a set map ${\P}^{top}(X)\longrightarrow
{\P}^{top}(Y)$. In other terms, one has the commutative diagram of
sets
\[\xymatrix{
{\P}^{top}(X)\fr{}\fd{\subset}& {\P}^{top}(Y)\fd{\subset}\\
\top([0,1],|X|) \fr{} &\top([0,1],|Y|).}
\]
\ed

\bd 
A {\rm globular complex } (resp. a {\rm relative globular complex})
$X$ is a globular precomplex (resp. a relative globular precomplex)
such that the attaching maps $\phi_\beta$ are non-decreasing. A
morphism of globular complexes is a morphism of globular precomplexes
which is non-decreasing. The category of globular complexes together
with the morphisms of globular complexes as defined above is denoted
by $\gltop$. The set $\gltop(X,Y)$ of morphisms of globular complexes
from $X$ to $Y$ equipped with the Kelleyfication of the compact-open
topology is denoted by $\tgltop(X,Y)$.
\ed

\bd 
Let $X$ be a globular complex. A point $\alpha$ of $X^0$ such that
there are no non-constant execution paths ending to $\alpha$
(resp. starting from $\alpha$) is called {\rm initial state} (resp.
{\rm final state}). More generally, a point of $X^0$ will be sometime
called {\rm a state} as well.
\ed

\bth (\cite{model2} Theorem~III.3.1) 
There exists a unique functor $cat:\gltop\longrightarrow\dtop$ such
that
\begin{enumerate}
\item if $X=X^0$ is a discrete globular complex, then $cat(X)$ is
the achronal flow $X^0$ (``achronal'' meaning with an empty path space), 
\item if $Z=\mathbf{S}^{n-1}$ or $Z=\mathbf{D}^{n}$ for some integer $n\geq 0$, 
then $cat(\glob^{top}(Z))=\glob(Z)$, 
\item for any globular complex $X$ with globular decomposition
$(n_\beta,\phi_\beta)_{\beta<\lambda}$, for any limit
ordinal $\beta\leq\lambda$, the canonical morphism of flows
\[\liminj_{\alpha<\beta} cat(X_\alpha)\longrightarrow
cat(X_\beta)\] is an isomorphism of flows, 
\item for any globular complex $X$ with globular decomposition
$(n_\beta,\phi_\beta)_{\beta<\lambda}$, for any
$\beta<\lambda$, one has the pushout of flows
\[\xymatrix{\glob(\mathbf{S}^{n_\beta-1})\fr{cat(\phi_\beta)}\fd{}& cat(X_\beta)\fd{}\\
\glob(\mathbf{D}^{n_\beta})\fr{} & cat(X_{\beta+1}).\cocartesien}\]
\end{enumerate}
\eth

The following theorem is important for the sequel:

\bth \label{remonte} The functor $cat$ induces a functor, still
denoted by $cat$ from $\gltop$ to $\hda\subset \dtop$ since its image
is contained in $\hda$. For any flow $X$ of $\hda$, there exists a
globular complex $Y$ such that $cat(U)=X$, which is constructed by
using the globular decomposition of $X$.  \eth

\bpf 
The construction of $U$ is made in the proof of \cite{model2}
Theorem~V.4.1.
\epf

\section{T-homotopy equivalence}
\label{rappelT}

The old notion of T-homotopy equivalence for globular complexes was
given in \cite{diCW}. A notion of T-homotopy equivalence of flows was
given in \cite{model2} and it was proved in the same paper that these
two notions are equivalent.

We first recall the definition of the branching and merging space
functors, and then the definition of a T-homotopy equivalence of
flows, exactly as given in \cite{model2} (Definition~\ref{Tdi} below),
and finally a characterization of T-homotopy of flows using globular
complexes (Theorem~\ref{TTT} below).

Roughly speaking, the branching space of a flow is the space of germs
of non-constant execution paths beginning in the same way.

\bp \label{universalpm} (\cite{exbranch} Proposition~3.1)
Let $X$ be a flow. There exists a topological space $\P^-X$ unique up
to homeomorphism and a continuous map $h^-:\P X\longrightarrow \P^- X$
satisfying the following universal property:
\begin{enumerate}
\item For any $x$ and $y$ in $\P X$ such that $t(x)=s(y)$, the equality
$h^-(x)=h^-(x*y)$ holds.
\item Let $\phi:\P X\longrightarrow Y$ be a
continuous map such that for any $x$ and $y$ of $\P X$ such that
$t(x)=s(y)$, the equality $\phi(x)=\phi(x*y)$ holds. Then there exists a
unique continuous map $\overline{\phi}:\P^-X\longrightarrow Y$ such that
$\phi=\overline{\phi}\circ h^-$.
\end{enumerate}
Moreover, one has the homeomorphism
\[\P^-X\iso \bigsqcup_{\alpha\in X^0} \P^-_\alpha X\]
where $\P^-_\alpha X:=h^-\lp \bigsqcup_{\beta\in
X^0} \P^-_{\alpha,\beta} X\rp$. The mapping $X\mapsto \P^-X$
yields a functor $\P^-$ from $\dtop$ to $\top$. 
\ep

\bd 
Let $X$ be a flow. The topological space $\P^-X$ is called the {\rm
branching space} of the flow $X$. The functor $\P^-$ is called the 
{\rm branching space functor}. \ed

\bp \label{universalmp}(\cite{exbranch} Proposition~A.1)
Let $X$ be a flow. There exists a topological space $\P^+X$ unique up
to homeomorphism and a continuous map $h^+:\P X\longrightarrow \P^+ X$
satisfying the following universal property:
\begin{enumerate}
\item For any $x$ and $y$ in $\P X$ such that $t(x)=s(y)$, the equality
$h^+(y)=h^+(x*y)$ holds.
\item Let $\phi:\P X\longrightarrow Y$ be a
continuous map such that for any $x$ and $y$ of $\P X$ such that
$t(x)=s(y)$, the equality $\phi(y)=\phi(x*y)$ holds. Then there exists a
unique continuous map $\overline{\phi}:\P^+X\longrightarrow Y$ such that
$\phi=\overline{\phi}\circ h^+$.
\end{enumerate}
Moreover, one has the homeomorphism
\[\P^+X\iso \bigsqcup_{\alpha\in X^0} \P^+_\alpha X\]
where $\P^+_\alpha X:=h^+\lp \bigsqcup_{\beta\in
X^0} \P^+_{\alpha,\beta} X\rp$. The mapping $X\mapsto \P^+X$
yields a functor $\P^+$ from $\dtop$ to $\top$. 
\ep

Roughly speaking, the merging space of a flow is the space of germs of
non-constant execution paths ending in the same way.

\bd \label{defplus}
Let $X$ be a flow. The topological space $\P^+X$ is called the {\rm
merging space} of the flow $X$. The functor $\P^+$ is called the {\rm
merging space functor}. \ed

\bd \cite{model2}
Let $X$ be a flow. Let $A$ and $B$ be two subsets of $X^0$. One says
that $A$ is {\rm surrounded} by $B$ (in $X$) if for any $\alpha\in A$,
either $\alpha \in B$ or there exists execution paths $\gamma_1$ and
$\gamma_2$ of $\P X$ such that $s(\gamma_1)\in B$,
$t(\gamma_1)=s(\gamma_2)=\alpha$ and $t(\gamma_2)\in B$. We denote
this situation by $A\lll B$.  \ed

\bd \cite{model2} 
Let $X$ be a flow. Let $A$ be a subset of $X^0$. Then the
\textit{restriction} $X\!\restriction_A$ of $X$ over $A$ is the unique
flow such that $(X\!\restriction_A)^0=A$, such that
$\P_{\alpha,\beta}(X\!\restriction_A)=\P_{\alpha,\beta}X$ for any
$(\alpha,\beta)\in A\p A$ and such that the inclusions $A\subset X^0$
and $\P (X\!\restriction_A)\subset \P X$ induces a morphism of flows
$X\!\restriction_A\longrightarrow X$.
\ed

\bd \label{Tdi}\cite{model2} 
Let $X$ and $Y$ be two objects of $\hda$.  A morphism of flows
$f:X\longrightarrow Y$ is a {\rm T-homotopy equivalence} if and only
if the following conditions are satisfied:
\begin{enumerate}
\item The morphism of flows $f:X\longrightarrow
Y\!\restriction_{f(X^0)}$ is an isomorphism of flows. In particular,
the set map $f^0:X^0\longrightarrow Y^0$ is one-to-one.
\item For $\alpha\in Y^0\backslash f(X^0)$, the topological spaces
$\P^-_\alpha Y$ and $\P^+_\alpha Y$ are singletons.
\item $Y^0\lll f(X^0)$.
\end{enumerate}
\ed

We recall the following important theorem for the sequel:

\bth (\cite{model2} Theorem~VI.3.5) \label{TTT} Let $X$ and $Y$ be two
objects of $\hda$. Let $U$ and $V$ be two globular complexes with
$cat(U)=X$ and $cat(V)=Y$ ($U$ and $V$ always exist by
Theorem~\ref{remonte}). Then a morphism of flows $f:X\longrightarrow Y$
is a T-homotopy equivalence if and only if there exists a morphism of
globular complexes $g:U\longrightarrow V$ such that $cat(g)=f$ and
such that the continuous map $|g|:|U|\longrightarrow |V|$ between the
underlying topological spaces is an homeomorphism.  \eth

This characterization was actually the first definition of a
T-homotopy equivalence proposed in \cite{diCW} (see Definition~4.10
p66).

\section{Some facts about relative $I^{gl}_+$-cell complexes}
\label{factI}

Recall that $I^{gl}_+=I^{gl}\cup
\{R:\{0,1\}\longrightarrow \{0\},C:\varnothing\longrightarrow \{0\}\}$
with
\[I^{gl}=\{\glob(\mathbf{S}^{n-1})\subset \glob(\mathbf{D}^{n}), n\geq
0\}.\] Let $I_g=I^{gl}\cup\{C\}$. Since for any $n\geq 0$, the
inclusion $\mathbf{S}^{n-1}\subset \mathbf{D}^n$ is a closed inclusion
of topological spaces, so an effective monomorphism of the category
$\top$ of compactly generated topological spaces, every morphism of
$I_g$, and therefore every morphism of $\cell(I_g)$, is an effective
monomorphism of flows as well (cf. also \cite{model3} Theorem~10.6).

\bp \label{removingR} 
If $f:X\longrightarrow Y$ is a relative $I^{gl}_+$-cell complex and if
$f$ induces a one-to-one set map from $X^0$ to $Y^0$, then
$f:X\longrightarrow Y$ is a relative $I_g$-cell subcomplex.
\ep

\bpf A pushout of $R$ appearing in the presentation of $f$ cannot
identify two elements of $X^0$ since, by hypothesis, $f^0:X^0
\rightarrow Y^0$ is one-to-one. So either such a pushout is trivial,
or it identifies two elements added by a pushout of $C$. 
\epf

\bp \label{factorR}
If $f:X\longrightarrow Y$ is a relative $I^{gl}_+$-cell complex, then
$f$ factors as a composite $g\circ h\circ k$ where $k:X\longrightarrow
Z$ is a morphism of $\cell(\{R\})$, where $h:Z \longrightarrow T$ is a
morphism of $\cell(\{C\})$, and where $g:T\longrightarrow Y$ is a
relative $I^{gl}$-cell complex.
\ep

\bpf 
One can use the small object argument with $\{R\}$ by \cite{model3}
Proposition~11.8. Therefore the morphism $f:X\longrightarrow Y$
factors as a composite $g\circ h$ where $h:X\longrightarrow Z$ is a
morphism of $\cell(\{R\})$ and where the morphism $Z\longrightarrow Y$
is a morphism of $\inj(\{R\})$. One deduces that the set map
$Z^0\longrightarrow Y^0$ is one-to-one. One has the pushout diagram of
flows
\[
\xymatrix{
X\fd{f} \fr{k} & Z \fd{g} \\
Y  \ar@{=}[r] & Y. \cocartesien}
\]
Therefore the morphism $Z\longrightarrow Y$ is a relative
$I^{gl}_+$-cell complex. Proposition~\ref{removingR} implies that the
morphism $Z\longrightarrow Y$ is a relative $I_g$-cell complex. The
morphism $Z\longrightarrow Y$ factors as a composite
$h:Z\longrightarrow Z \sqcup (Y^0\backslash Z^0)$ and the inclusion
$g:Z \sqcup (Y^0\backslash Z^0)\longrightarrow Y$. 
\epf

\bp\label{tjrcof} 
Let $X=X^0$ be a set viewed as a flow (i.e. with an empty path space).
Let $Y$ be an object of $\hda$. Then any morphism from $X$ to $Y$ is a
cofibration.
\ep

\bpf 
Let $f:X\longrightarrow Y$ be a morphism of flows. Then $f$ factors as
a composite $X=X^0\longrightarrow Y^0\longrightarrow Y$.  Any set map
$X^0\longrightarrow Y^0$ is a transfinite composition of pushouts of
$C$ and $R$. So any set morphism $X^0\longrightarrow Y^0$ is a
cofibration of flows. And for any flow $Y$, the canonical morphism of
flows $Y^0\longrightarrow Y$ is a cofibration since it is a relative
$I_g$-cell complex. Hence the result.
\epf

\section{Left properness of the weak S-homotopy model structure of $\dtop$}
\label{sectionleftproperflow}

\bp\label{pushexplicit}\label{model3} (\cite{model3} Proposition~15.1)
Let $f:U\longrightarrow V$ be a continuous map.  Consider the pushout
diagram of flows
\[
\xymatrix{
\glob(U)\fr{}\fd{\glob(f)} & X \fd{g}\\
\glob(V)\fr{} & Y. \cocartesien}
\]
Then the continuous map $\P g:\P X\longrightarrow \P Y$ is a
transfinite composition of pushouts of continuous maps of the form a
finite product $\id\p\dots\p f\p\dots \p\id$ where the symbol $\id$
denotes identity maps.
\ep

\bp \label{versleft}
Let $f:U\longrightarrow V$ be a Serre cofibration. Then the pushout of
a weak homotopy equivalence along a map of the form a finite product
$\id_{X_1}\p\dots\p f\p\dots \p\id_{X_p}$ with $p\geq 0$ is still a
weak homotopy equivalence.
\ep

If the topological spaces $X_i$ for $1\leq i\leq p$ are cofibrant,
then the continuous map $\id_{X_1}\p\dots\p f\p\dots \p\dots\p
\id_{X_p}$ is a cofibration since the model category of compactly
generated topological spaces is monoidal with the categorical product
as monoidal structure.  So in this case, the result follows from the
left properness of this model category (\cite{ref_model2}
Theorem~13.1.10). In the general case, $\id_{X_1} \p \dots \p f\p\dots
\p \dots \p \id_{X_p}$ is not a cofibration anymore. But any
cofibration $f$ for the Quillen model structure of $\top$ is a
cofibration for the Str{\o}m model structure of $\top$
\cite{MR35:2284} \cite{MR39:4846} \cite{ruse} \cite{strom2}. In the
latter model structure, any space is cofibrant. Therefore the
continuous map $\id_{X_1} \p \dots \p f \p \dots \p \dots \p
\id_{X_p}$ is a cofibration of the Str{\o}m model structure of $\top$,
that is a NDR pair. So the continuous map $\id_{X_1} \p \dots \p f \p
\dots \p \dots \p \id_{X_p}$ is a closed $T_1$-inclusion anyway. This
fact will be used below.

\bpf 
We already know that the pushout of a weak homotopy equivalence along
a cofibration is a weak homotopy equivalence. The proof of this
proposition is actually an adaptation of the proof of the left
properness of the model category of compactly generated topological
spaces. Any cofibration is a retract of a transfinite composition of
pushouts of inclusions of the form $\mathbf{S}^{n-1}\subset
\mathbf{D}^{n}$ for $n\geq 0$. Since the category of compactly generated
topological spaces is cartesian closed, the binary product preserves
colimits. Thus we are reduced to considering a diagram of topological
spaces like
\[
\xymatrix{
X_1\p \dots \p \mathbf{S}^{n-1}\p\dots \p X_p\fd{} \fr{} & U \fd{} \fr{s}& X \fd{} \\
 X_1\p \dots \p \mathbf{D}^{n} \p\dots\p X_p\fr{} & \widehat{U} \cocartesien \fr{\widehat{s}} & \widehat{X} \cocartesien}
\]
where $s$ is a weak homotopy equivalence and we have to prove that
$\widehat{s}$ is a weak homotopy equivalence as well. By
\cite{MR36:6480} and \cite{MR48:6288}, it suffices to prove that
$\widehat{s}$ induces a bijection between the path-connected
components of $\widehat{U}$ and $\widehat{X}$, a bijection between the
fundamental groupoids $\pi(\widehat{U})$ and $\pi(\widehat{X})$, and
that for any local coefficient system of abelian groups $A$ of
$\widehat{X}$, one has the isomorphism
$\widehat{s}^*:H^*(\widehat{X},A)\iso H^*(\widehat{U},\widehat{s}^*A)$.

For $n=0$, one has $\mathbf{S}^{n-1}=\varnothing$ and
$\mathbf{D}^{n}=\{0\}$.  So $X_1\p \dots \p \mathbf{S}^{n-1}\p\dots \p
X_p = \varnothing$ and $X_1\p \dots \p \mathbf{D}^{n} \p \dots \p X_p
= X_1\p \dots \p X_p$. So $\widehat{U}\iso U\sqcup (X_1\p \dots \p
X_p)$ and $\widehat{X}\iso X\sqcup (X_1\p \dots \p X_p)$.  Therefore
the mapping $t$ is the disjoint sum $s\sqcup \id_{X_1\p \dots \p
  X_p}$. So it is a weak homotopy equivalence.

Let $n\geq 1$. The assertion concerning the path-connected components
is clear. Let $\mathbf{T}^n = \{x\in \mathbb{R}^n,0<|x|\leq 1\}$.
Consider the diagram of topological spaces
\[
\xymatrix{
X_1\p \dots \p \mathbf{S}^{n-1}\p\dots \p X_p\fd{} \fr{} & U \fd{} \fr{s}& X \fd{} \\
 X_1\p \dots \p \mathbf{T}^{n} \p \dots \p X_p\fr{} & \widetilde{U} \cocartesien \fr{\widetilde{s}} & \widetilde{X}. \cocartesien}
\]
Since the pair $(\mathbf{T}^{n},\mathbf{S}^{n-1})$ is a deformation
retract, the three pairs $(X_1\p \dots \p \mathbf{T}^{n} \p \dots \p
X_p,X_1\p \dots
\p \mathbf{S}^{n-1}\p\dots \p X_p)$, $(\widetilde{U},U)$ and
$(\widetilde{X},X)$ are deformation retracts as well. So the
continuous maps $U\longrightarrow \widetilde{U}$ and $X\longrightarrow
\widetilde{X}$ are both homotopy equivalences.
The Seifert-Van-Kampen theorem for the fundamental groupoid
(cf. \cite{MR48:6288} again) then yields the diagram of groupoids
\[
\xymatrix{
\pi(X_1\p \dots \p \mathbf{T}^{n}\p\dots \p X_p)\fd{} \fr{} & \pi(\widetilde{U}) \fd{} \fr{\pi(\widetilde{s})}& \pi(\widetilde{X}) \fd{} \\
 \pi(X_1\p \dots \p \mathbf{D}^{n} \p X_p)\fr{} & \pi(\widehat{U})
 \cocartesien \fr{\pi(\widehat{s})} & \pi(\widehat{X}). \cocartesien}
\]
Since $\pi(\widetilde{s})$ is an isomorphism of groupoids, then so is $\pi(\widehat{s})$.

Let $\mathbf{B}^n=\{x\in
\mathbb{R}^n,0\leq |x|< 1\}$.  Then $(\mathbf{B}^n,\widetilde{U})$ is an
excisive pair of $\widehat{U}$ and $(\mathbf{B}^n,\widetilde{X})$ is an
excisive pair of $\widehat{X}$. The Mayer-Vietoris long exact sequence then yields
the commutative diagram of groups
\[
\xymatrix{
\dots\fr{}& H^p(\widehat{X},A)\fr{} \fd{} & H^p({X},A)\oplus H^p(\mathbf{B}^n,A)\fr{}\fd{\iso} &
H^p(\mathbf{B}^n\backslash\{0\},A)\fr{}\fd{\iso}&\dots\\
\dots\fr{}& H^p(\widehat{U},\widehat{s}^*A)\fr{} & H^p({U},s^*A)\oplus H^p(\mathbf{B}^n,s^*A)\fr{} &
H^p(\mathbf{B}^n\backslash\{0\},s^*A)\fr{}&\dots}
\]
A five-lemma argument completes the proof.
\epf

\bp \label{petitplus} Let $\lambda$ be an ordinal. Let $M:\lambda
\longrightarrow \top$ and $N:\lambda \longrightarrow \top$ be two
$\lambda$-sequences of topological spaces. Let $s:M\longrightarrow N$
be a morphism of $\lambda$-sequences which is also an objectwise weak
homotopy equivalence. Finally, let us suppose that for all
$\mu<\lambda$, the continuous maps $M_\mu \longrightarrow M_{\mu+1}$
and $N_\mu \longrightarrow N_{\mu+1}$ are of the form a finite product
$\id_{X_1}\p\dots\p f\p\dots \p\id_{X_p}$ with $p\geq 0$ and with $f$
a Serre cofibration.  Then the continuous map $\liminj s:\liminj M
\longrightarrow \liminj N$ is a weak homotopy equivalence.  \ep

If for all $\mu<\lambda$, the continuous maps $M_\mu \longrightarrow
M_{\mu+1}$ and $N_\mu \longrightarrow N_{\mu+1}$ are cofibrations,
then Proposition~\ref{petitplus} above is a consequence of
\cite{ref_model2} Proposition~17.9.3 and of the fact that the model
category $\top$ is left proper. With the same additional hypotheses,
Proposition~\ref{petitplus} above is also a consequence of
\cite{MR2045835} Theorem~A.7. Indeed, the latter states that an
homotopy colimit can be calculated either in the usual Quillen model
structure of $\top$, or in the Str{\o}m model structure of $\top$
\cite{ruse} \cite{strom2}.

\bpf The principle of the proof is standard.  If the ordinal $\lambda$
is not a limit ordinal, then this is a consequence of
Proposition~\ref{versleft}. Assume now that $\lambda$ is a limit
ordinal. Then $\lambda \geq \aleph_0$.

Let $u:\mathbf{S}^n\longrightarrow \liminj N$ be a continuous map.
Then $u$ factors as a composite $\mathbf{S}^n\longrightarrow N_\mu
\longrightarrow \liminj N$ since the $n$-dimensional sphere
$\mathbf{S}^n$ is compact and since any compact space is
$\aleph_0$-small relative to closed $T_1$-inclusions
(\cite{MR99h:55031} Proposition~2.4.2). By hypothesis, there exists a
continuous map $\mathbf{S}^n \longrightarrow M_\mu$ such that the
composite $\mathbf{S}^n \longrightarrow M_\mu \longrightarrow N_\mu$
is homotopic to $\mathbf{S}^n\longrightarrow N_\mu$. Hence the
surjectivity of the set map $\pi_n(\liminj M,*) \longrightarrow
\pi_n(\liminj N,*)$ (where $\pi_n$ denotes the $n$-th homotopy group)
for $n\geq 0$ and for any base point $*$.

Let $u,v:\mathbf{S}^n \longrightarrow \liminj M$ be two continuous
maps such that there exists an homotopy $H:\mathbf{S}^n\p [0,1]
\longrightarrow \liminj N$ between $\liminj s \circ f$ and $\liminj s
\circ g$. Since the space $\mathbf{S}^n\p [0,1]$ is compact, the homotopy 
$H$ factors as a composite $\mathbf{S}^n\p [0,1] \longrightarrow
N_{\mu_0} \longrightarrow \liminj N$ for some $\mu_0<\lambda$. And
again since the space $\mathbf{S}^n$ is compact, the map $f$
(resp. $g$) factors as a composite $\mathbf{S}^n \longrightarrow
M_{\mu_1} \longrightarrow \liminj M$ (resp. $\mathbf{S}^n
\longrightarrow M_{\mu_2}
\longrightarrow \liminj M$) with $\mu_1<\lambda$ (resp. $\mu_2<\lambda$). 
Then $\mu_4=\max(\mu_0,\mu_1,\mu_2)<\lambda$ since $\lambda$ is a
limit ordinal. And the map $H:\mathbf{S}^n\p [0,1] \longrightarrow
N_{\mu_4}$ is an homotopy between $f:\mathbf{S}^n \longrightarrow
M_{\mu_4}$ and $g:\mathbf{S}^n \longrightarrow M_{\mu_4}$. So the set
map $\pi_n(\liminj M,*) \longrightarrow \pi_n(\liminj N,*)$ for $n\geq
0$ and for any base point $*$ is one-to-one. 
\epf

\bth \label{leftproperflow} 
The model category $\dtop$ is left proper.
\eth

\bpf Consider the pushout diagram of $\dtop$
\[
\xymatrix{
U \fr{s}\fd{i} & X \fd{} \\ V \fr{t} & \cocartesien Y}
\]
where $i$ is a cofibration of $\dtop$ and $s$ a weak S-homotopy
equivalence. We have to check that $t$ is a weak S-homotopy
equivalence as well. The morphism $i$ is a retract of a
$I^{gl}_+$-cell complex $j:U\longrightarrow W$. If one considers
the pushout diagram of $\dtop$
\[\xymatrix{
U \fr{s}\fd{j} & X \fd{} \\
W \fr{u} & \cocartesien Y}
\]
then $t$ must be a retract of $u$. Therefore it suffices to prove that
$u$ is a weak S-homotopy equivalence. So one can suppose that one has
a diagram of flows of the form
\[
\xymatrix{
A\fr{\phi} \fd{k} & U \fr{s}\fd{i} & X \fd{} \\
B \fr{} & \cocartesien V \fr{t} & \cocartesien Y}
\]
where $k\in \cell(I^{gl}_+)$. By Proposition~\ref{factorR}, the
morphism $k:A\longrightarrow B$ factors as a composite $A
\longrightarrow A' \longrightarrow A'' \longrightarrow B$ where the morphism 
$A \longrightarrow A'$ is an element of $\cell(\{R\})$, where the
morphism $A' \longrightarrow A''$ is an element of $\cell(\{C\})$, and
where the morphism $A'' \longrightarrow B$ is a morphism of
$\cell(I^{gl})$. So we have to treat the cases $k\in \cell(\{R\})$,
$k\in \cell(\{C\})$ and $k\in \cell(I^{gl})$.

The case $k\in \cell(I^{gl})$ is a consequence of
Proposition~\ref{pushexplicit}, Proposition~\ref{versleft} and
Proposition~\ref{petitplus}. The case $k\in \cell(\{C\})$ is trivial.

Let $k\in \cell(\{R\})$. Let $(\alpha,\beta)\in U^0\p U^0$. Then 
$\P_{i(\alpha),i(\beta)}V$ (resp.
$\P_{i(\alpha),i(\beta)}Y$) is a coproduct of terms of the form
$\P_{\alpha,u_0}U\p \P_{v_0,u_1}U\p \dots\p \P_{v_p,\beta}U$
(resp. $\P_{\alpha,u_0}X\p \P_{v_0,u_1}X\p \dots\p
\P_{v_p,\beta}X$) such that $(u_i,v_i)$ is a pair of distinct elements of 
$U^0=X^0$ identified by $k$. So $t$ is a weak S-homotopy equivalence
since a binary product of weak homotopy equivalences is a weak
homotopy equivalence.
\epf

\section{T-homotopy equivalence and $I^{gl}_+$-cell complex}
\label{main1}

The first step to understand the reason why Definition~\ref{Tdi} is
badly-behaved is the following theorem which gives a description of
the T-homotopy equivalences $f:X\longrightarrow Y$ such that the
$0$-skeleton of $Y$ contains exactly one more state than the
$0$-skeleton of $X$.

\bth \label{re1} 
Let $X$ and $Y$ be two objects of $\hda$. Let $f:X\longrightarrow Y$ be
a T-homotopy equivalence. Assume that $Y^0=X^0\sqcup \{\alpha\}$. Then
the canonical morphism $\varnothing\longrightarrow X$ factors as a
composite $\varnothing\longrightarrow u_f(X)\longrightarrow
v_f(X)\longrightarrow X$ such that
\begin{enumerate}
\item one has the diagram 
\[
\xymatrix{
& \varnothing \fd{}&\\
\{\widehat{0},\widehat{1}\}=\glob(\mathbf{S}^{-1}) \fr{} \fd{} & u_f(X) \fd{}&&\\
\vI=\glob(\mathbf{D}^0) \fr{}\fd{\phi} & \cocartesien v_f(X)  \fd{} \ar@{->}[rr] && X \fd{}\\
\vI*\vI \fr{} & \cocartesien \widehat{v}_f(X) \ar@{->}[rr] && \cocartesien Y
}
\]
\item the morphisms $\varnothing\longrightarrow u_f(X)$ and 
$v_f(X)\longrightarrow X$ are relative $I_g$-cell complexes. 
\end{enumerate}
\eth

By Proposition~\ref{tjrcof}, the morphism
$\{\widehat{0},\widehat{1}\}=\glob(\mathbf{S}^{-1})\longrightarrow
u_f(X)$ is a cofibration.  Therefore the morphism $\vI\longrightarrow
v_f(X)$ is a cofibration as well. The morphism $u_f(X)\longrightarrow
\widehat{v}_f(X)$ is a relative $I_g$-cell complex as well since it is
a pushout of the inclusion $\{\widehat{0},\widehat{1}\}\subset
\vI*\vI$ sending $\widehat{0}$ to the initial state of $\vI*\vI$ and $\widehat{1}$ to the
final state of $\vI*\vI$.

\bpf By Proposition~\ref{removingR}, and since $Y$ is an object of
$\hda$, the canonical morphism of flows $Y^0\longrightarrow Y$ is a
relative $I_g$-cell complex. So there exists an ordinal $\lambda$ and
a $\lambda$-sequence $\mu \mapsto Y_\mu :\lambda\longrightarrow \dtop$
(so also denoted by $Y$) such that $Y=\liminj_{\mu<\lambda} Y_\mu$ and
such that for any ordinal $\mu<\lambda$, the morphism
$Y_\mu\longrightarrow Y_{\mu+1}$ is a pushout of the form
\[
\xymatrix{
\glob(\mathbf{S}^{n_\mu-1})\fd{} \fr{\phi_\mu}& Y_\mu\fd{}\\
\glob(\mathbf{D}^{n_\mu}) \ar@{->}[r]_-{\psi_\mu}& Y_{\mu+1}\cocartesien
}
\]
of the inclusion of flows $\glob(\mathbf{S}^{n_\mu})\longrightarrow
\glob(\mathbf{D}^{n_\mu+1})$ for some $n_\mu\geq 0$.

For any ordinal $\mu$, the morphism of flows $Y_\mu\longrightarrow
Y_{\mu+1}$ induces an isomorphism between the $0$-skeletons $Y_\mu^0$
and $Y_{\mu+1}^0$. If $n_\mu\geq 1$ for some $\mu$, then for any
$\beta,\gamma\in Y^0$, the topological space $\P_{\beta,\gamma}Y_\mu$
is non empty if and only if the topological space
$\P_{\beta,\gamma}Y_{\mu+1}$ is non empty. Consider the set of
ordinals
\[\left\{\mu<\lambda; \bigsqcup_{\beta\in X^0}\P_{\beta,\alpha}Y_\mu\neq \varnothing\right\}\]
It is non-empty since $f$ is a T-homotopy equivalence. Take its smallest element
$\mu_0$. Consider the set of ordinals
\[\left\{\mu<\lambda; \bigsqcup_{\beta\in X^0}\P_{\alpha,\beta}Y_\mu\neq \varnothing\right\}\]
Take its smallest element $\mu_1$. Let us suppose for instance that
$\mu_0<\mu_1$.

The ordinal $\mu_0$ cannot be a limit ordinal. Otherwise for any
$\mu<\mu_0$, the isomorphisms of flows $Y_\mu=Z_\mu\sqcup \{\alpha\}$
and $Y_{\mu_0}\iso \liminj_{\mu<\mu_0}\lp Z_\mu\sqcup
\{\alpha\}\rp\iso\lp
\liminj_{\mu<\mu_0}Z_\mu\rp \sqcup \{\alpha\}$ would hold:
contradiction. Therefore $\mu_0=\mu_2+1$ and $n_{\mu_2}=0$. There
does not exist other ordinal $\mu$ such that
$\phi_\mu(\widehat{1})=\alpha$ otherwise $\P_\alpha^+Y$ could not be a
singleton anymore.

For a slightly different reason, the ordinal $\mu_1$ cannot be a limit
ordinal either. Otherwise if $\mu_1$ was a limit ordinal, then the
isomorphism of flows $Y_{\mu_1}\iso \liminj_{\mu<\mu_1} Y_\mu$ would
hold. The path space of a colimit of flows is in general not the
colimit of the path spaces. But any element of $\P Y_{\mu_1}$ is a
composite $\gamma_1*\dots*\gamma_p$ where the $\gamma_i$ for $1\leq
i\leq p$ belong to $\liminj_{\mu<\mu_1}\P Y_\mu$. By hypothesis, there
exists an execution path $\gamma_1*\dots*\gamma_p\in
\P_{\alpha,\beta}Y_{\mu_1}$ for some $\beta\in X^0$.  So
$s(\gamma_1)=\alpha$, which contradicts the definition of
$\mu_1$. Therefore $\mu_1=\mu_3+1$ and necessarily $n_{\mu_3}=0$.
There does not exist any other ordinal $\mu$ such that
$\phi_\mu(\widehat{0})=\alpha$ otherwise $\P_\alpha^-Y$ could not be a
singleton anymore.

Therefore one has the following situation: $Y_{\mu_2}$ is a flow of
the form $Z_{\mu_2}\sqcup \{\alpha\}$. The passage from $Y_{\mu_2}$ to
$Y_{\mu_2+1}$ is as follows:
\[
\xymatrix{
\glob(\mathbf{S}^{-1}) \fd{} \fr{\phi_{\mu_2}} & Y_{\mu_2} \fd{}\\
\glob(\mathbf{D}^0) \fr{} & \cocartesien Y_{\mu_2+1}}
\]
where $\phi_{\mu_2}(\widehat{0})\in X^0$ and
$\phi_{\mu_2}(\widehat{1})=\alpha$. The morphism of flows
$Y_{\mu_2+1}\longrightarrow Y_{\mu_3}$ is a transfinite composition of
pushouts of the inclusion of flows
$\glob(\mathbf{S}^{n})\longrightarrow \glob(\mathbf{D}^{n+1})$ where
$\phi_\mu(\widehat{0})$ and $\phi_\mu(\widehat{1})$ are never equal to
$\alpha$. The passage from $Y_{\mu_3}$ to $Y_{\mu_3+1}$ is as follows:
\[
\xymatrix{
\glob(\mathbf{S}^{-1}) \fd{} \fr{\phi_{\mu_3}} & Y_{\mu_3} \fd{}\\
\glob(\mathbf{D}^0) \fr{} & \cocartesien Y_{\mu_3+1}}
\]
where $\phi_{\mu_3}(\widehat{0})=\alpha$ and
$\phi_{\mu_3}(\widehat{1})\in X^0$.  The morphism of flows
$Y_{\mu_3+1}\longrightarrow Y_{\lambda}$ is a transfinite composition
of pushouts of the inclusion of flows
$\glob(\mathbf{S}^{n})\longrightarrow \glob(\mathbf{D}^{n+1})$ where
$\phi_\mu(\widehat{0})$ and $\phi_\mu(\widehat{1})$ are never equal to
$\alpha$.  Hence the result.
\epf

We are now ready to give a characterization of the old T-homotopy
equivalences:

\bth 
\label{re3}
Let $X$ and $Y$ be two objects of $\hda$.  Then a morphism of flows
$f:X\longrightarrow Y$ is a T-homotopy equivalence if and only if
there exists a commutative diagram of flows of the form (with $\vI^{*
(n+1)}:=\vI^{* n}*\vI$ and $\vI^{* 1}:=\vI$ for $n\geq 1$):
\[
\xymatrix{
& \varnothing \fd{}&\\
\bigsqcup_{i\in I}\{0,1\}=\bigsqcup_{i\in I}\glob(\mathbf{S}^{-1}) \fr{} \fd{} & u_f(X) \fd{}&&\\
\bigsqcup_{i\in I}\vI=\bigsqcup_{i\in I}\glob(\mathbf{D}^0) \fr{}\fd{\bigsqcup_{i\in I} r_i} & \cocartesien v_f(X)  \fd{} \ar@{->}[rr] && X \fd{f}\\
\bigsqcup_{i\in I}\vI^{* n_i} \fr{} & \cocartesien \widehat{v}_f(X) \ar@{->}[rr] && \cocartesien Y
}
\]
where for any $i\in I$, $n_i$ is an integer with $n_i\geq 1$ and such
that $r_i:\vI\longrightarrow \vI^{* n_i}$ is the unique morphism of
flows preserving the initial and final states and where 
the morphisms $\varnothing\longrightarrow u_f(X)$ and 
$v_f(X)\longrightarrow X$ are relative $I_g$-cell complexes. 
\eth

The pushout above tells us that the copy of $\vI$ corresponding to the
indexing $i\in I$ is divided in the concatenation of $n_i$ copies of
$\vI$. This intuitively corresponds to a refinement of observation.

\bpf
By Theorem~\ref{remonte}, there exists a globular complex $U$
(resp. $V$) such that $cat(U)=X$ (resp. $cat(V)=Y$). If a morphism of
flows $f:X\longrightarrow Y$ is a T-homotopy equivalence, then by
Theorem~\ref{TTT}, there exists a morphism of globular complexes
$g:U\longrightarrow V$ such that $cat(g)=f$ and such that the
continuous map $|g|:|U|\longrightarrow |V|$ between the underlying
topological spaces is an homeomorphism. So for any pair of points
$(\alpha,\beta)$ of $X^0\p X^0$, and any morphism $\vI\longrightarrow
X$ appearing in the globular decomposition of $X$, the set of
subdivision of this segment in $Y$ is finite since $Y^0$ is discrete
and since the segment $[0,1]$ is compact. The result is then
established by repeatedly applying Theorem~\ref{re1}.

Now suppose that a morphism of flows $f:X=cat(U)\longrightarrow
Y=cat(V)$ can be written as a pushout of the form of the statement of
the theorem. Then start from a globular decomposition of $U$ which is
compatible with the composite $\varnothing\longrightarrow
u_f(X)\longrightarrow v_f(X)\longrightarrow X$. Then let us divide
each segment of $[0,1]$ corresponding to the copy of $\vI$ indexed by
$i\in I$ in $n_i$ pieces. Then one obtains a globular decomposition of
$V$ and the identity of $U$ gives rise to a morphism of globular
complexes $g:U\longrightarrow V$ which induces an homeomorphism
between the underlying topological spaces and such that
$cat(g)=f$. Hence the result.
\epf

\bd Let $n\geq 1$. The {\rm full $n$-cube} $\vec{C}_n$ is by
definition the flow $Q(\{\widehat{0}<\widehat{1}\}^n)$, where $Q$ is
the cofibrant replacement functor. \ed

The flow $\vec{C}_3$ is represented in Figure~\ref{3cube}.

\begin{lem} \label{ordrestate} If a flow $X$ is loopless, then the
  transitive closure of the set \[\{(\alpha,\beta)\in X^0\p X^0\hbox{
    such that }\P_{\alpha,\beta}X\neq\varnothing\}\] induces a partial
  ordering on $X^0$.
\end{lem}

\bpf If $(\alpha,\beta)$ and $(\beta,\alpha)$ with $\alpha\neq \beta$
belong to the transitive closure, then there exists a finite sequence
$(x_1,\dots,x_\ell)$ of elements of $X^0$ with $x_1=\alpha$,
$x_\ell=\alpha$, $\ell>1$ and with $\P_{x_m,x_{m+1}}X$ non-empty for
each $m$. Consequently, the space $\P_{\alpha,\alpha}X$ is non-empty
because of the existence of the composition law of $X$: contradiction.
\epf

\bth \label{pasiso} 
Let $n\geq 3$. There does not exist any zig-zag sequence
\[
\xymatrix{
\vec{C}_n=X_0 \fr{f_0} & X_1 & X_2 \fl{f_1}\fr{f_2} & \dots &  X_{2n}=\vI \fl{f_{2n-1}}
}
\]
where every $X_i$ is an object of $\hda$ and where every $f_i$ is
either a S-homotopy equivalence or a T-homotopy equivalence. \eth

\bpf By an immediate induction, one sees that each flow $X_i$ is loopless, 
with a finite $0$-skeleton. Moreover by construction, each poset
$(X_i^0,\leq)$ is bounded, i.e. with one bottom element $\widehat{0}$
and one top element $\widehat{1}$. So the zig-zag sequence above gives
rise to a zig-zag sequence of posets:
\[
\xymatrix{
\vec{C}_n^0=X_0^0=\{\widehat{0}<\widehat{1}\}^n \fr{} & X_1^0 &
X_2^0 \fl{}\fr{} & \dots & X_{2n}^0=\vI^0=\{\widehat{0}<\widehat{1}\}
\fl{} }
\]
where $\{\widehat{0}<\widehat{1}\}^n$ is the product
$\{\widehat{0}<\widehat{1}\}\p \dots \p \{\widehat{0}<\widehat{1}\}$
($n$ times) in the category of posets.  Each morphism of posets is an
isomorphism if the corresponding morphism of flows is a S-homotopy
equivalence because a S-homotopy equivalence induces a bijection
between the $0$-skeletons. Otherwise one can suppose by
Theorem~\ref{re1} that the morphism of posets $P_1\longrightarrow P_2$
can be described as follows: take a segment $[x,y]$ of $P_1$ such that
$]x,y[=\varnothing$; add a vertex $z\in ]x,y[$; then let
$P_2=P_1\cup\{z\}$ with the partial ordering $x<z<y$.  In such a
situation, $\min(]z,-[)$ exists and is equal to $y$, and $\max(]-,z[)$
exists and is equal to $x$. So by an immediate induction, there must
exist $x,y,z\in \{0<1\}^n$ with $x<z<y$ and such that $\min(]z,-[)=y$
and $\max(]-,z[)=x$. This situation is impossible in the poset
$\{\widehat{0}<\widehat{1}\}^n$ for $n\geq 3$.
\epf

\section{Generalized T-homotopy equivalence}
\label{main2}

As explained in the introduction, it is not satisfactory not to be
able to identify $\vec{C}_3$, and more generally $\vec{C}_n$ for
$n\geq 3$, with $\vI$. The following definitions are going to be
important for the sequel of the paper, and also for the whole series.

\bd~\footnote{The statement of the definition is slightly different, 
but equivalent to the statement given in other parts of this series.}
A {\rm full directed ball} is a flow $\vec{D}$ such that:
\begin{itemize}
\item $\vec{D}$ is loopless (so by Lemma~\ref{ordrestate}, the set
  $\vec{D}^0$ is equipped with a partial ordering $\leq$)
\item $(\vec{D}^0,\leq)$ is finite bounded
\item for all $(\alpha,\beta)\in \vec{D}^0\p \vec{D}^0$, the
  topological space $\P_{\alpha,\beta}\vec{D}$ is weakly contractible
  if $\alpha<\beta$, and empty otherwise by definition of $\leq$.
\end{itemize}
\ed

Let $\vec{D}$ be a full directed ball. Then by Lemma~\ref{ordrestate},
the set $\vec{D}^0$ can be viewed as a finite bounded poset.
Conversely, if $P$ is a finite bounded poset, let us consider the
\textit{flow} $F(P)$ \textit{associated with} $P$: it is of course
defined as the unique flow $F(P)$ such that $F(P)^0=P$ and
$\P_{\alpha,\beta}F(P)=\{u\}$ if $\alpha<\beta$ and
$\P_{\alpha,\beta}F(P)=\varnothing$ otherwise. Then $F(P)$ is a full
directed ball and for any full directed ball $\vec{D}$, the two flows
$\vec{D}$ and $F(\vec{D}^0)$ are weakly S-homotopy equivalent.

Let $\vec{E}$ be another full directed ball. Let $f: \vec{D}
\longrightarrow \vec{E}$ be a morphism of flows preserving the initial
and final states. Then $f$ induces a morphism of posets from
$\vec{D}^0$ to $\vec{E}^0$ such that $f(\min
\vec{D}^0) = \min \vec{E}^0$ and $f(\max \vec{D}^0) = \max \vec{E}^0$. 
Hence the following definition:

\bd \label{definitiondeT}
Let $\mathcal{T}$ be the class of morphisms of posets
$f:P_1\longrightarrow P_2$ such that:
\begin{enumerate}
\item The posets $P_1$ and $P_2$ are finite and bounded. 
\item The morphism of posets $f:P_1 \longrightarrow P_2$ is one-to-one; 
in particular, if $x$ and $y$ are two elements of $P_1$ with $x<y$,
then $f(x)<f(y)$.
\item One has $f(\min P_1)=\min P_2$ and  $f(\max P_1)=\max P_2$.
\end{enumerate}
Then a {\rm generalized T-homotopy equivalence} is a morphism of
$\cof(\{Q(F(f)),f\in\mathcal{T}\})$ where $Q$ is the cofibrant
replacement functor of the model category $\dtop$.
\ed

It is of course possible to identity $\vec{C}_n$ ($n\geq 1$) with
$\vI$ by the following zig-zag sequence of S-homotopy and generalized
T-homotopy equivalences 
\[
\xymatrix{ \vI & \fl{\simeq} Q(\vI) \ar@{->}[rr]^-{Q(F(g_n))} &&
  Q(\{\widehat{0} < \widehat{1}\}^n), }
\] 
where $g_n:\{\widehat{0} < \widehat{1}\} \longrightarrow \{\widehat{0} <
\widehat{1}\}^n \in \mathcal{T}$.

The relationship between the new definition of T-homotopy equivalence
and the old definition is as follows:

\bth  \label{OKnew}
Let $X$ and $Y$ be two objects of $\hda$.  Let $f:X\longrightarrow Y$
be a T-homotopy equivalence. Then $f$ can be written as a composite
$X\longrightarrow Z\longrightarrow Y$ where $g:X\longrightarrow Z$ is
a generalized T-homotopy equivalence and where $h:Z\longrightarrow Y$
is a weak S-homotopy equivalence. 
\eth

\bpf By Theorem~\ref{re3}, there exists a pushout
diagram of flows of the form (with $\vI^{* (n+1)}:=\vI^{* n}*\vI$ and
$\vI^{* 1}:=\vI$ for $n\geq 1$):
\[
\xymatrix{
\bigsqcup_{k\in K} \vI \ar@{^{(}->}[r]^{}\fd{\bigsqcup_{k\in K} r_k} & X \fd{}\\
\bigsqcup_{k\in K} \vI^{* n_k} \fr{} & Y \cocartesien}
\]
where for any $k\in K$, $n_k$ is an integer with $n_k\geq 1$ and such
that $r_k:\vI\longrightarrow \vI^{* n_k}$ is the unique morphism of
flows preserving the initial and final states. Notice that each
$\vI^{* n_k}$ is a full directed ball. Thus one obtains the following 
commutative diagram: 
\[
\xymatrix{
\bigsqcup_{k\in K} Q(\vI) \ar@{->>}[r]^{\simeq}\fd{} & \bigsqcup_{k\in K} \vI 
\ar@{^{(}->}[r]\fd{} & X \ar@{^{(}->}[d]\\ 
\bigsqcup_{k\in K} Q(\vI^{* n_k}) \ar@{=}[r] & \bigsqcup_{k\in K} Q(\vI^{* n_k}) \ar@{^{(}->}[r]\fd{\simeq} & \cocartesien Z \fd{\simeq} \\
&\bigsqcup_{k\in K} \vI^{* n_k} \ar@{^{(}->}[r] & Y. \cocartesien}
\] 
Now here are some justifications for this diagram. First of all, a
morphism of flows $f:M\longrightarrow N$ is a fibration of flows if
and only if the continuous map $\P f :\P M\longrightarrow \P N$ is a
Serre fibration of topological spaces. Since any coproduct of Serre
fibration is a Serre fibration, the morphism of flows $\bigsqcup_{i\in
I} Q(\vI)\longrightarrow \bigsqcup_{k\in K} \vI $ is a trivial
fibration of flows. Thus the underlying set map $\bigsqcup_{k\in K}
Q(\vI)\longrightarrow \bigsqcup_{k\in K} \vI $ is surjective. So the
commutative square
\[
\xymatrix{
\bigsqcup_{k\in K} Q(\vI)\ar@{->}[d]\fr{} & \fd{} X\\
\bigsqcup_{k\in K} Q(\vI^{* n_k})\fr{} & \cocartesien Z}
\]
is cocartesian and the morphism of flows $X\longrightarrow Z$ is then
a generalized T-homotopy equivalence. It is clear that the morphism
$\bigsqcup_{k\in K} Q(\vI^{* n_k})\longrightarrow
\bigsqcup_{k\in K} \vI^{* n_k}$ is a weak S-homotopy equivalence. The latter 
morphism is even a fibration of flows, but that does not matter
here. So the morphism $Z\longrightarrow Y$ is the pushout of a weak
S-homotopy equivalence along the cofibration $\bigsqcup_{k\in K}
Q(\vI^{* n_k})\longrightarrow Z$.  Since the model category $\dtop$ is
left proper by Theorem~\ref{leftproperflow}, the proof is complete.
\epf

\section{Conclusion}

This new definition of T-homotopy equivalence contains the old one up
to S-homotopy equivalence.  The drawback of the old definition
presented in \cite{model2} is overcome. It is proved in \cite{3eme}
that this new notion of T-homotopy equivalence does preserve the
branching and merging homology theories. And it is proved in
\cite{4eme} that the underlying homotopy type of a flow is also
preserved by this new definition of T-homotopy equivalence. Finally,
the paper \cite{hocont} proposes an application of this new notion of
dihomotopy: that is a Whitehead theorem for the full dihomotopy
relation.

\end{document}